\documentclass[reqno]{amsart}

\newtheorem{theorem}{Theorem}[section]
\newtheorem{lemma}[theorem]{Lemma}

\theoremstyle{definition}
\newtheorem{definition}[theorem]{Definition}
\newtheorem{proposition}[theorem]{Proposition}

\newtheorem{corollary}[theorem]{Corollary}

\theoremstyle{remark}
\newtheorem{remark}[theorem]{Remark}

\newcommand{\id}{\mbox{id}} 
 
\newcommand{\Ker}{\mbox{Ker\,}} 

\newcommand{\End}{\mbox{End}} 
\newcommand{\Hom}{\mbox{Hom}}

\newcommand{\eps}{\varepsilon}

\newcommand{\actl}{\leftharpoonup}
\newcommand{\lact}{\triangleright} 
 
\newcommand{\la}{\langle\,} 
\newcommand{\ra}{\,\rangle}

\newcommand{\0}{_{(0)}}  
\newcommand{\1}{_{(1)}} 
\newcommand{\2}{_{(2)}} 
\newcommand{\3}{_{(3)}}

\title[Hopf algebra actions on Separable
Extensions of Depth 2]{Hopf Algebra Actions on Strongly 
Separable Extensions of Depth Two}

\author{Lars Kadison}
\address{Chalmers University of Technology/G{\" o}teborg
University \\ Matematiskt Centrum \\
S-412 96 G{\" o}teborg \\
Sweden}
\email{lkadison@online.no}
\thanks{The first author thanks A.A.~Stolin for useful discussions,
and NorFA in Oslo for financial support.}

\author{Dmitri Nikshych}
\address{U.C.L.A., Department of Mathematics, Los Angeles, CA 90095-1555, USA}
\email{nikshych@math.ucla.edu}
\date{}

\begin{document}
\begin{abstract}
We  bring together ideas in analysis of Hopf
$*$-algebra actions on II$_1$ subfactors of finite Jones index
\cite{GHJ, S} and algebraic characterizations of Frobenius, Galois  
and cleft Hopf extensions \cite{KT,K,DT} to prove a
non-commutative algebraic analogue of the classical theorem:
a finite field extension is Galois iff it is separable and normal.
Suppose $N \hookrightarrow M$ is a separable Frobenius extension
of $k$-algebras split as $N$-bimodules with a trivial centralizer $C_M(N)$.
Let $M_1 := \End(M_N)$ and $M_2 := \End(M_1)_M$
be the endomorphism algebras in the Jones tower $N \hookrightarrow M
\hookrightarrow M_1 \hookrightarrow M_2$.  
We show that under depth $2$
conditions on the second centralizers $A := C_{M_1}(N)$
and $B : = C_{M_2}(M)$ the algebras $A$ and $B$ are semisimple
Hopf algebras dual to one another and such that $M_1$ is a smash product
of $M$ and $A$, and  that $M$ is a $B$-Galois extension of $N$.

\end{abstract}

\maketitle

\begin{section}
{Introduction}
Three well-known functors associated to the induced representations
of a subalgebra pair $N \subseteq M$ are restriction $\mathcal{R}$ of $M$-modules
to $N$-modules, its adjoint $\mathcal{T}$ which tensors $N$-modules by $M$,
and its co-adjoint $\mathcal{H}$ which applies $\Hom_N(M,-)$ to $N$-modules. 
The algebra extension $M/N$ is said to be \textit{Frobenius} if
$\mathcal{T}$ is naturally isomorphic to $\mathcal{H}$ \cite{Mo}.  $M/N$ is said to be
\textit{separable} if the counit of adjunction $\mathcal{T}\mathcal{R} \stackrel{\cdot}{\rightarrow} 1$
 is  naturally split epi;
and $M/N$ is a \textit{split extension} if the unit of adjunction $1 
\stackrel{\cdot}{\rightarrow}
\mathcal{R}\mathcal{T}$ is  naturally split
monic \cite{P}. An algebraic model for finite
Jones index subfactor theory is given in  \cite{K1,K2} using a \textit{strongly separable extension}, 
which has  all three
of these properties. Over a ground field $k$,  an
 \textit{irreducible} extension $M/N$, which is characterized by having
trivial centralizer $C_M(N) = k1$,
is  strongly separable if it is   split, separable and Frobenius.

In this paper we extend the results of Szyma\'{n}ski \cite{S}
and others \cite{Lo,EN} on Hopf $*$-algebra 
actions and finite index subfactors with  a trace (i.e., linear map
$\phi: M \rightarrow k$ such that $\phi(mm') = \phi(m'm)$ for all 
$m, m' \in M$ and $\phi(1) = 1_k$) 
 to  strongly separable, irreducible extensions.  We will
not require of our algebras that they possess a  trace.
However, we require some hypotheses on the endomorphism algebra
$M_1$ of the natural module $M_N$, to which there is a monomorphism
given by the left regular representation of $M$ in $\End(M_N)$.  
We require a \textit{depth two condition} that two successive
endomorphism algebra extensions, $M \hookrightarrow M_1$
and $M_1 \hookrightarrow M_2$, be free with bases in the second
centralizers, $A := C_{M_1}(N)$ and $B := C_{M_2}(M)$. 
Working over a field of arbitrary
characteristic, we prove in Theorem~\ref{smash product}
and Theorem~\ref{M_1 is a smash product}: 
\begin{theorem}
\label{summa summarum}
The Jones tower $M \subseteq M_1 \subseteq M_2$ over a strongly separable, 
irreducible   extension $N \subseteq M$ of  depth 2  has centralizers $A$
and $B$ that are involutive semisimple Hopf algebras dual to 
one another, with an action of $B$ on $M_1$
and another action of $A$ on $M$ such that $M_1$ and $M_2$ are smash products:
$M_2 \cong M_1 \# B$ and $M_1 \cong M \# A$.
\end{theorem}
 
The main theorem~\ref{summa summarum} is of intrinsic interest in 
 extending \cite{S} to the case of an
irreducible finite index  pair of von Neumann factors of arbitrary
type (I, II, or III). 
%%%new
Secondly, it  gives a proof that $M_1$ is
a smash product without appeal to a tunnel construction; i.e. assuming
the strong hypotheses in 
a characterization of a strongly separable 
extension $M/N$ that is the endomorphism 
algebra extension of some $N/R$. 
Thirdly, the main theorem is the difficult 
piece in the proof of a non-commutative analogue of the classical theorem 
in field theory \cite{St}:

\begin{theorem}
\label{classical theorem}
A finite field extension $E'/F'$ is Galois if and only if 
$E'/F'$ is separable and normal.
\end{theorem}

By $E'/F'$  Galois we mean that the Galois group 
$G$ of $F'$-algebra automorphisms of $E'$ has $F'$ as its
fixed field ${E'}^G$. From a modern point of view, the right non-commutative generalization 
of Galois extension is the Hopf-Galois extension (cf.\ Section~3) \cite{M},
with classical Galois groups interpreted as 
cosemisimple Hopf algebras. From the modern
 cohomological point of view, the non-commutative
 separable extensions mentioned above are a direct
generalization of separable field extensions (cf.\ Section 2) \cite{HS}.
The trace map $T: E' \rightarrow F'$
for finite separable field extensions \cite{L} is a Frobenius homomorphism for a Frobenius
extension (cf.\ Section~2), while the trace map for Galois extensions
is the action of an integral on the overfield
(corresponding to the mapping $E$ in the proof of Theorem~\ref{converse}). 
In Section~6 we prove
 the following non-commutative analogue of Theorem~\ref{classical theorem}:

\begin{theorem}
\label{non-commutative analogue}
If $M/N$ is an irreducible  extension of depth 2, then $M/N$
is  strongly separable if and only if $M/N$ is an $H$-Galois extension, where
$H$ is a semisimple, cosemisimple Hopf algebra. 
\end{theorem} 
 
 In Sections 2 and 3 we note that the non-commutative notions of separable extension
and Hopf-Galois extension  generalize 
separability and Galois extension, respectively, for finite field extensions. 
However, Theorem~\ref{non-commutative analogue} 
is not a generalization of the classical theorem, since non-trivial field extensions
are not irreducible. The proof of 
Theorem~\ref{non-commutative analogue}  follows from
Theorems~\ref{converse} and~\ref{summa summarum}.  Theorem~\ref{converse} is 
an easier result with roots in 
\cite{KT,D,K2}.  The smash
product result for $M_2$ in Theorem~\ref{summa summarum}
follows from the depth 2 properties in Section~3, the non-degenerate pairing
of $A$ and $B$ in Section~4, and the action of $B$ on $M_1$ in Section~5
together with the key Proposition~\ref{exchange relation}.  The non-degenerate pairing
in Eq.\ (\ref{pairing}) transfers the algebra structures of $A$ and $B$ 
to coalgebra structures on $B$ and $A$, respectively, that result in the Hopf algebra
structures on these. The antipodes on $A$ and $B$
come from a basic symmetry in the definition of the pairing. 
From the action of $B$ on $M_1$ with fixed subalgebra $M$, 
we dualize in Section~6 to an  $A$-extension $M_1/M$, compute that
it is $A$-cleft, and use
the Hopf algebra-theoretic characterization of the latter 
as a crossed product : we show that $M_1$ is a smash
product of $M$ with $A$ from the triviality of the cocycle.
Each section begins with an introduction
to the main  terminology, theory and results in the section.    
 
\end{section}

\begin{section}
{Strongly separable extensions with trivial centralizer}

In this section, we recall the most basic definitions and facts
for irreducible and split extensions, Frobenius extensions and algebras,  separable extensions and algebras,
and strongly separable extensions and algebras.  We introduce Frobenius
homomorphisms and their dual bases, which characterize Frobenius extensions,
noting that Frobenius homomorphisms are faithful, and have Nakayama automorphisms
measuring their deviation from being a trace on the centralizer. 
After introducing separability and strongly separable extensions, we
come to the  important theory of the basic
construction $M_1$, conditional expectation $E_M: M_1 \rightarrow M$
and Jones idempotent $e_1 \in M_1$. The basic construction is repeated
to form the tower of algebras $N \subseteq M  \subseteq M_1  \subseteq M_2$,
and the braid-like relations between $e_1$ and $e_2 \in M_2$ are pointed
out.

Throughout this paper, $k$ denotes a field.  Let $M$ and $N$ be  associative unital
 $k$-algebras with $N$ a unital subalgebra of $M$.   We refer to
  $N \subseteq M$ or a (unity-preserving) monomorphism $N \hookrightarrow M$
 as an {\it algebra extension} $M/N$.  We note
the endomorphism algebra extension $\End(M_N)/M$ obtained
from $m \rightarrow \lambda_m$ for each $m \in M$, where $\lambda_m$
is  left multiplication by $m \in M$, a right $N$-module endomorphism of $M$.

In this section, we denote the \textit{centralizer} of a bimodule ${}_NP_N$ by 
$P^N := \{ p \in P | \, \forall n \in N, \, pn = np \}$, a special
case of which is the centralizer subalgebra of $N$ in $M$:  $C_M(N) = M^N$.
The algebra extension $M/N$ will be called {\it irreducible} 
if the centralizer subalgebra 
 is trivial, i.e.,  $C_M(N) = k1$.  
In this case the centers $Z(M)$ and $Z(N)$ both lie in $C_M(N)$, so they are  
trivial as well. If $\mathcal{E}$ denotes $\End(M_N)$ and
 $M^{op}$ denotes the opposite algebra of  $M$, we note that $(\forall m \in M$)
\begin{equation}
C_{\mathcal{E}}(M) = \{ f \in \mathcal{E}|\, mf(x) = f(mx), \forall m \in M \} = 
\End({}_M M_N) \cong C_M(N)^{op}.
\label{endo}
\end{equation}
Whence the endomorphism algebra extension is irreducible too. 

$M/N$ is a \textit{split extension} if there is an $N$-bimodule projection
$E: M \rightarrow N$.  Thus, $E(1) = 1$, $E(nmn') = nE(m)n'$,
for all $n,n' \in N, m \in M$, and $M = N \oplus \ker E$ as
$N$-bimodules, the last being an equivalent condition.  The condition
mentioned in the first paragraph of Section~1 is easily shown to be equivalent as well \cite{P}. 

\subsection*{Frobenius extensions}
$M/N$ is said to be a { \it Frobenius extension} if the natural right
$N$-module $M_N$ is finitely generated projective and there is a 
bimodule isomorphism of $M$ with
its (algebra extension) dual:  ${}_N M_M \cong {}_N\Hom(M_N,N_N)_M$ \cite{K}. 
This definition 
is equivalent to the condition that $M/N$ has
 a bimodule homomorphism $E: {}_N M_N 
\rightarrow {}_N N_N$, called a {\it Frobenius homomorphism}, and  
elements in $M$, 
$\{ x_i \}_{i=1}^n$, 
$\{ y_i \}_{i=1}^n$, called {\it dual bases}, such that the equations 
\begin{equation}
\sum_{i=1}^n E(mx_i)y_i = m  = \sum_{i=1}^n x_i E(y_im)
\label{Frobeq}
\end{equation}
hold for every $m \in M$ \cite{K}.\footnote{For if $\{x_i\}$, $\{ f_i\}$ is a projective base
for $M_N$ and $E$ is the image of $1$, then there is $y_i \mapsto Ey_i = f_i$
such that $\sum_i x_iEy_i = \id_M$. The other equation follows.  Conversely, $M_N$ is 
explicitly finitely generated projective, while 
$x \mapsto Ex$ is bijective.} In particular,  Frobenius extension
may be defined equivalently in terms of the natural {\it left} module ${}_N M$ instead.
The Hattori-Stallings rank of the projective modules $M_N$ or ${}_N M$ are
both given by $\sum_i E(y_ix_i)$ in $N/[N,N]$ \cite{K1}. 
It is not hard to check that the \textit{index} 
$[M:N]_E := \sum_i x_i y_i \in Z(M)$ 
(use Eqs.\ (\ref{Frobeq}))  depends only on $E$, and $E(1) \in Z(N)$. 
Furthermore, $M/N$ is split if and only if there is a $d \in C_M(N)$ such
that $E(d) = 1$ \cite{K1}.  

If $M_N$ is  free,  $M/N$ is called a 
{\em free Frobenius extension} \cite{K}.  
By choosing dual bases $\{ x_i \}$, $\{ f_i \}$ for $M_N$
such that $f_i(x_j) = \delta_{ij}$, 
we arrive at {\it orthogonal dual bases}
$\{x_i \}$, $\{ y_i \}$, which satisfy 
$
E(y_i x_j) = \delta_{ij} 
$.
Conversely, with $E$, $x_i$ and $y_i$ satisfying this equation, it is clear that
$M/N$ is free Frobenius. 

If $N$ is the unit subalgebra $k1$, $M$ is a \textit{Frobenius algebra}, 
a notion introduced  in a 1903 paper of Frobenius \cite{F}.  
Such an algebra $M$ is characterized 
by having a \textit{faithful}, or non-degenerate,
 linear functional $E: M \rightarrow k$; i.e., $E(Mm) = 0$ implies 
$m = 0$, or equivalently, $E(mM) = 0$ implies $m = 0$
(in one direction a trivial
application of Eqs.\ (\ref{Frobeq})). 

We note the following \textit{transitivity} result
 with an easy proof. Consider the tower of algebras 
$N \subseteq M \subseteq R$.  If $M/N$ and $R/M$ are  Frobenius extensions,
then so is the composite extension $R/N$.  Moreover, the following
proposition has a proof left to the reader:  

\begin{proposition}
\label{transitivity}
If $M/N$ and $R/M$ are algebra extensions with Frobenius homomorphisms $E: M \rightarrow N$,
$F: R \rightarrow M$ and
dual bases $\{ x_i \}$, $\{ y_i \}$  and
$\{ z_j \}$, $\{w_j \}$, respectively, then 
$R/N$  has Frobenius homomorphism $E \circ F$
and dual bases $\{ z_j x_i \}$, $\{ y_i w_j \} $.  
\end{proposition}

If $M/N$ and $R/M$ are irreducible, the composite index satisfies
the Lagrange equation:
\[ [R: N]_{EF} = [R:M]_F [M:N]_E.
\] 

\subsection*{Nakayama automorphism}
Given a Frobenius homomorphism $E : M \rightarrow N$ and
 an element $c$ in the centralizer $C_M(N)$, the maps
$cE$ and $Ec$ defined by $cE(x) := E(xc)$ and $Ec(x) = E(cx)$
are both $N$-bimodule maps belonging to 
the $N$-centralizers of both the $N$-bimodules $\Hom_N(M_N,N_N)$
and $\Hom_N({}_NM,{}_NN)$.  Since $m \mapsto Em$ is a bimodule isomorphism, 
 ${}_NM_M \cong {}_N\Hom_N^r(M,N)_M$, it follows that
there is a unique $c' \in C_M(N) = M^N$ such that $Ec' = cE$.  The mapping $q: c \mapsto c'$
on $C_M(N)$ is clearly an automorphism, 
called the \textit{Nakayama automorphism}, or
modular automorphism, with defining equation given by 
\begin{equation}
E(q(c)m) = E(mc)
\label{Naka}
\end{equation}
for every $c \in C_M(N)$ and $m \in M$ \cite{K}. 
$M/N$ is a \textit{symmetric} Frobenius extension if $q$ is an inner automorphism.
In case $N = k1$, this recovers the usual notion of \textit{symmetric algebra} (a finite-dimensional
algebra with non-degenerate or faithful trace), for 
if $q: M \rightarrow M$ is given by $q(m) = umu^{-1}$, 
then $Eu$ is such a trace by Eq.\ (\ref{Naka}).

\subsection*{Separability}
Throughout this paper we consider $ M \otimes_N M $
with its natural $M$-$M$-bimodule structure.
$M/N$ is said to be a {\it separable extension} if the multiplication
epimorphism $\mu: M \otimes_N M \rightarrow M$ has a right inverse  as  $M$-$M$-bimodule
homomorphisms \cite{HS}. This is clearly equivalent to the existence of an element
$e \in M \otimes_N M$ such that $me = em$ for every $m \in M$
and $\mu(e) = 1$, called a {\it separability
element}: separable extensions are precisely the algebra extensions with trivial 
relative Hochschild cohomology groups in degree one or more \cite{HS}. 
A Frobenius extension $M/N$ with $E,x_i,y_i$ as before is separable
if and only if there is a $d \in C_M(N)$ such that $\sum_i x_i d y_i = 1$
\cite{HS}. 

If $N = k1_M$, $M/N$ is a separable extension iff $M$ is 
a separable $k$-algebra; i.e. a finite dimensional, semisimple $k$-algebra with matrix
blocks over division algebras $D_i$ where  $Z(D_i)$ is a finite separable (field)
 extension of $k$.
If $k$ is algebraically closed, each $D_i = k$ and $M$ is isomorphic to a direct product of
matrix blocks of order $n_i$ 
over $k$. 

For example, if $E'/F'$ is a finite separable field extension, $\alpha \in E'$
the primitive element such that $E' = F'(\alpha)$, and $p(x) =  x^n - \sum_{i=0}^{n-1}
c_i x^i$ the minimal polynomial of $\alpha$ in $F'[X]$, then a separability
element is given by  
\[   \sum_{i=0}^{n-1} \alpha^{i}\otimes_{F'}\frac{\sum_{j=0}^{i} 
c_{j}\alpha^{j}}{p'(\alpha)\alpha^{i+1}}. 
\]  

A $k$-algebra $M$ is said to be {\it strongly separable} in Kanzaki's sense
 if $M$ has a {\it symmetric} separability
element $e$ (necessarily unique); i.e., $\tau(e) = e$ where $\tau$ is the twist map on $M \otimes_k M$.
An equivalent condition  is that $M$ has a trace $t : M \rightarrow k$ (i.e.,
$t(mn) = t(nm)$ for all $m, n \in M$) and elements $x_1,\ldots,x_n,y_1,\ldots,y_n$
such that $\sum_i t(mx_i)y_i = m$ for all $m \in M$ 
and $\sum_i x_i y_i = 1_M$. A third equivalent
condition is that $M$ has an invertible Hattori-Stalling rank over its center \cite{De}. It
follows that 
the characteristic of $k$ does not divide the orders $n_i$ of the matrix
blocks (i.e., $n_i 1_k \neq 0$); for a separable $k$-algebra $M$,
this is also a sufficient condition for strong separability
 in case $k$ is algebraically closed.

\subsection*{Strongly separable extensions} We are now ready to 
define the main object of investigation in this paper.

\begin{definition}[cf.\ \cite{K1,K2}] 
\label{strongly separable extension}
A $k$-algebra extension $N \subseteq M$ is called a {\em strongly 
separable, irreducible 
extension}
 if $M/N$ is an irreducible Frobenius extension with Frobenius homomorphism
 $E: M \rightarrow N$, and dual bases $\{ x_i \}$, $\{ y_i \}$
 such that 
\begin{enumerate}
\item $E(1) \neq 0$,
\item $\sum_i x_i y_i \neq 0 $,
\end{enumerate} 

\end{definition}

\begin{remark}
Since $M/N$ is irreducible, the centers of $M$ and $N$
are trivial, so $E(1) = \mu 1_S$ for some nonzero $\mu \in k$. 
Then $\frac{1}{\mu} E, \mu x_i, y_i$ is a new Frobenius homomorphism
with dual bases for $M/N$. With no loss of generality then,
we assume that
\begin{equation}
E(1) = 1.
\label{normalization}
\end{equation}
It follows that $M = N \oplus \Ker E$ as $N$-$N$-bimodules and $E^2 = E$
when $E$ is viewed in $\End_N(M)$. 
Also 
\begin{equation}
\sum_i x_i y_i = \lambda^{-1} 1_M
\label{index}
\end{equation}
for some nonzero $\lambda \in k$.  It follows that
$\lambda \sum_i x_i \otimes y_i$ is a separability element and $M/N$
is separable.  The data $E, x_i, y_i $ for a strongly separable, irreducible extension, 
satisfying Eqs.\ (\ref{normalization}) and (\ref{Frobeq}), 
is uniquely determined.\footnote{There is a  close but
complicated relationship between   Kanzaki strongly separable $k$-algebras
and strongly separable extensions $A/k1$ in the sense of
 \cite{K2}.  Note that $A = M_2(F_2)$,
where $F_2$ is a field of characteristic $2$, is not Kanzaki
strongly separable, but is a strongly separable extension $A/F_21$ since  
$E(A) = a_{11} + a_{12} + a_{21}$  and $$ \sum_i x_i \otimes y_i = 
e_{11} \otimes e_{21} + e_{12} \otimes e_{11} + e_{12} \otimes e_{21} + e_{22} \otimes e_{12}
+ e_{22} \otimes e_{22} + e_{21} \otimes e_{22},$$ satisfies
 $\sum_i x_i y_i = 1$, $E(1) = 1$, 
 $E$  a Frobenius homomorphism with dual bases $x_i, y_i$.  
However, a strongly
separable extension $A/k1$ with Markov trace \cite{K2} is  Kanzaki 
strongly separable; and conversely, if $k = Z(A)$.}
\end{remark}

\subsection*{The basic construction.} The basic construction
begins with the following \textit{endomorphism ring theorem},
whose proof we sketch here for the sake of completeness: 
\begin{theorem}[Cf.\ \cite{K1,K2}]
\label{endomorphism ring theorem}
$\mathcal{E}/M$ is a strongly separable, irreducible
extension of index $\lambda^{-1}$.
\end{theorem}
\begin{proof}
For a Frobenius extension $M/N$, we have $\mathcal{E} \cong M \otimes_N M$
by sending $ f \mapsto \sum_i f(x_i) \otimes y_i$ with inverse
$m \otimes n \mapsto \lambda_m E \lambda_n$ in the notation above. 
We denote $M_1 := M \otimes_N M$, and note that
the multiplication on  $M_1$ induced by composition
of endomorphisms is given by the \textit{$E$-multiplication}:
\begin{equation}
(m_1 \otimes m_2)(m_3 \otimes m_4) = m_1 E(m_2 m_3) \otimes m_4.
\label{E-multiplication}
\end{equation}
The unity element is $1_1 := \sum_i x_i \otimes y_i$ 
in the notation above. It is easy to see that $E_M := \lambda \mu$,
where $\mu$ is the multiplication mapping $M_1 \rightarrow M$,
is a normalized Frobenius homomorphism, and $\{ \lambda^{-1} x_i \otimes 1 \}$,
$\{ 1 \otimes y_i \}$ are  dual bases satisfying
Eqs.\ (\ref{normalization}) and (\ref{index}). 
\end{proof}

We make note of the {\it first Jones idempotent}, $e_1 := 1 \otimes 1 \in M_1$,
which cyclically generates $M_1$ as an $M$-$M$-bimodule: $M_1 = \{ \sum_i x_i e_1 y_i |
\, x_i, y_i \in M \}$. In this paper, a Frobenius homomorphism
$E$ satisfying $E(1) = 1$ is called a {\it conditional expectation}. 
We describe $M_1, e_1,  E_M$ as the ``basic construction''
of $N \subseteq M$.

\subsection*{The Jones tower.}  The basic construction is repeated in order to produce
the Jones tower of $k$-algebras above $N \subseteq M$:
\begin{equation}
N \subseteq M \subseteq M_1 \subseteq M_2 \subseteq \cdots
\label{Jones tower}
\end{equation}
In this paper
we will only need to consider $M_2$, which is the basic construction of 
$M \subseteq M_1$.
As such it is given by 
\begin{equation}
M_2 = M_1 \otimes_M M_1 \cong M \otimes_N M \otimes_N M 
\label{emtwo}
\end{equation}
 with $E_M$-multiplication, and conditional expectation 
$E_{M_1} := \lambda \mu: M_2 \rightarrow M_1$ given
by 
\[
m_1 \otimes m_2 \otimes m_3 \mapsto \lambda m_1 E(m_2) \otimes m_3.
\]
The second Jones idempotent is given by
\[
e_2 = 1_1 \otimes 1_1 = \sum_{i,j} x_i \otimes y_i x_j \otimes y_j,
\]
and satisfies $e_2^2 = e_2$ in the $E_M$-multiplication of $M_2$. 

\subsection*{The braid-like relations.}  Note that $1_2 = \sum_i \lambda^{-1}
x_i \otimes 1 \otimes y_i$ and
\(
E_{M_i}(e_{i+1}) = \lambda 1 \)
where $M_0 $ denotes $M$. Then the following relations between $e_1, e_2$ are 
readily computed in $M_2$ without the hypothesis of irreducibility: 
\begin{proposition}
\begin{eqnarray*}
e_1 e_2 e_1 & = & \lambda e_1  1_2 \\
e_2 e_1 e_2 & = & \lambda e_2.
\end{eqnarray*}
\end{proposition}
\begin{proof}  The proof can be found in \cite[Ch.~3]{K1}.
\end{proof}

%%%%%%%%%%%%%%%%%%%%%%%%%%%%%%%%%%%%%%%%%%%%%%%%%%%%%%%%%%%%%%%%%%%%%%%%%%%%
\section{Depth 2 properties}

In this section, we place depth 2 conditions on the modules ${}_M M_1$
and ${}_{M_1}M_2$ by requiring that they are  free with bases in
$A := C_{M_1}(N)$ and $B := C_{M_2}(M)$, respectively. 
We then show that $A$ and $B$ are separable algebras with $E_M|_A$
and $E_{M_1}|_B$, respectively, as faithful linear functionals.  
The classical depth 2 property, coming from subfactor theory \cite{GHJ},
is established for the large centralizer, $C := C_{M_2}(N)$; i.e.,
$C$ is the basic construction of $A$ or $B$ over the trivial centralizer
with conditional expectations $E_A$ and $E_B$ studied later in the section.
We next establish the important property that $F := E_M \circ E_{M_1}$
restricts to a faithful linear functional on $C$. We interpret
the various Nakayama automorphisms arising from $F$, $E_M|_A$
and $E_{M_1}|_B$. The important Pimsner-Popa identities are
established. We end
this section by recalling the basic properties of Hopf-Galois extensions,
and prove Theorem~\ref{converse} which states that an $H$-Galois extension is
strongly separable of depth 2 if $H$ is a semisimple, cosemisimple Hopf
algebra. This establishes one of the implications in Theorem~\ref{non-commutative
analogue}.

\subsection*{Finite depth and  depth 2 conditions.}   
We extend the notion of  {\em depth} known in subfactor
theory \cite{GHJ} to Frobenius extensions. 

\begin{lemma}
\label{dual bases in A}
For all $n\geq 1$ in the Jones tower~(\ref{Jones tower})
the following conditions are equivalent (we denote $M_{-1} =N$
and $M_0= M$) :
\begin{enumerate}
\item[(1)] 
$M_{n-1}$ is a free right
$M_{n-2}$-module with a basis in $C_{M_{n-1}}(N)$
(respectively, $M_{n}$ is a free right $M_{n-1}$-module 
with a basis in $C_{M_n}(M)$).
\item[(2)] 
There exist orthogonal dual bases for $E_{M_{n-2}}$
in $C_{M_{n-1}}(N)$ (respectively,
there exist orthogonal dual bases for $E_{M_{n-1}}$ in $C_{M_n}(M)$).
\end{enumerate}
\end{lemma} 
\begin{proof}
We show that   (1) implies (2), the other implication is trivial.
Denote by $\{ z_i \}$ and $\{ w_i \}$  orthogonal dual bases in $M_{n-1}$ for 
$E_{M_{n-2}}$, where $\{ z_i \} \subset C_{M_{n-1}}(N)$.  
We compute that $w_i \in C_{M_{n-1}}(N)$:
$$
xw_i = \sum_j xE_{M_{n-2}}(w_iz_j)w_j = \sum_j \delta_{ij}xw_j 
= \sum_j E_{M_{n-2}}(w_i xz_j)w_j = w_ix
$$
for every $x \in N$.  The second statement in the proposition is 
proven similarly with dual bases $\{ u_j \}$ in $C_{M_n}(M)$ and therefore 
$\{ v_j \}$ in $C_{M_n}(M)$. 
%We fix the notation in this proof.
\end{proof}

We say that a Frobenius extension $M/N$ has
a {\em finite depth} if the equivalent conditions
of Lemma~\ref{dual bases in A} are satisfied for some $n\geq 1$ .
It is not hard to check that  in this case
they also hold true for $n+1$ (and, hence, for all
$k\geq n$). Indeed, if $\{ u_j \}$ and $\{ v_j \}$ are as above, then
$\{ \lambda^{-1} u_j e_{n+1}\},\,\{ e_{n+1}v_j \} \subset C_{M_{n+1}}(M)$
is a pair of orthogonal dual bases for $E_{M_n}$. 
We then  define the {\em depth} of a finite depth extension
$M/N$ to be the smallest number $n$ for which these conditions hold.
In the trivial case, an irreducible extension of depth $1$ leads
to $M = N$.

Let $A$ and $B$ denote the ``second'' centralizer algebras:
\[
A := C_{M_1}(N), \ \ \ \ B := C_{M_2}(M).
\]
The  {\em depth 2 conditions}
that we will use in this paper are then explicitly:
\begin{enumerate}
\item  $M_1$ is a free right $M$-module with  basis in $A$;
\item $M_2$ is a free right $M_1$-module with basis in $B$.
\end{enumerate}
It is easy to show that $M_1$ and $M_2$ are also free as left
$M-$ and $M_1-$modules, respectively. Note that
the depth 2 conditions make sense for an arbitrary ring extension
$M/N$ where $M_1$ and $M_2$ stand for the successive endomorphism rings. 

In what follows, we assume that $M/N$ has depth $2$ and denote 
$\{ z_i \},\, \{ w_i \} \subset A$   orthogonal dual bases for 
$E_M$ and $\{ u_i \},\, \{ v_i \} \subset B$
orthogonal dual bases for $E_{M_1}$ that exist by 
Lemma~\ref{dual bases in A}. 

\begin{proposition}
\label{A is  separable}
$A$ and $B$ are  separable algebras.
\end{proposition}
\begin{proof}
For all $a \in A$, we have 
$ \sum_i E_M(az_i)w_i = a = \sum_i z_i E_M(w_i a)
$
where $E_M(az_i)$ and $ E_M(w_i a)$ lie in $C_M(N) = k1_M$.  $\{ z_i \}$ is linearly independent
over $M$, whence over $k$, so $A$, similarly $B$, is  finite dimensional.

It follows that $E_M$ restricted to $A$ is a Frobenius homomorphism.  Since $\{ z_i \}$, $\{w_i \}$ are
dual bases and $[M_1:M]_{E_M} = \lambda^{-1}$, it follows  that $\lambda \sum_i z_i \otimes w_i$
is a separability element. 
Similarly, $B$ is a Frobenius algebra
with Frobenius homomorphism $E_{M_1}$, and a separable algebra 
with separability element $\lambda \sum_j u_j \otimes v_j$.  
\end{proof}
The lemma below is a first step to the main result 
that  
$M_2$ is a smash product of $B$ and $M_1$ 
(cf.\ Theorem~\ref{smash product}). 
\begin{lemma}
\label{tensor product}
We have 
$
  M_1 \cong M \otimes_k A
$
as $M$-$A$-bimodules, and
$
M_2 \cong M_1 \otimes_k B
$
as $M_1$-$B$-bimodules. 
\end{lemma}
\begin{proof}
We map $w \in M_1$ into $\sum_i E_M(wz_i) \otimes w_i \in M \otimes A $, which
has inverse mapping $m \otimes a \in M \otimes A$ into $ma \in M_1$. 

The proof of the second statement is completely similar. 
\end{proof}

We let $C = C_{M_2}(N)$.  Note that $A \subseteq C$ and $B \subseteq C$.
Of course $A1_2 \cap B = k1_2$ since $C_{M_1}(M) = k1_1$.  We
will now show in a series of steps  the classical depth 2 property 
that $C$ is the basic construction of $A$ or $B$
over the trivial centralizer. 

\begin{lemma}
\label{C is A tensor B}
$C \cong A \otimes_k B$  via multiplication 
$a \otimes b \mapsto ab$ and $C \cong B \otimes_k A$
via $b \otimes a \mapsto ba$.
%%%new  We need both isomorphisms in 4.6.
\end{lemma}
\begin{proof}
If $c \in C$, then $\sum_j E_{M_1}(cu_j) \otimes v_j \in A \otimes B$, which
provides an inverse to the first map above. The second part is established
similarly.
\end{proof}

\begin{lemma}
\label{e_2A is e_2C}
We have $e_2 A = e_2 C$ and $Ae_2 = Ce_2$ as subsets of $M_2$. 
Also, $e_1B = e_1 C$ and $B e_1 = C e_1$ in $M_2$. 
\end{lemma}
\begin{proof}
For each $b \in B$ we have $b_j,b'_j \in M_1$ such that 
$$
e_2 b =  1_1 \otimes 1_1 \sum_j b_j \otimes b'_j 
       =  e_2 \sum_j E_M(b_j)b'_j  \in ke_2
$$
since $\sum_j E_M(b_j)b'_j \in C_{M_1}(M) = k1$. 
Then
\(
e_2 C = e_2BA = e_2 A. \) The second equality is proven similarly.
The second statement is proven in the same way by making use of
$e_1 A =Ae_1 =  ke_1$. 
\end{proof}

We place the $E_M$-multiplication on $A \otimes A$,
and the $E_{M_1}$-multiplication on $B \otimes B$ below. 
 \begin{proposition}[Depth 2 property]
\label{C is A tensor A}
We have $C = Ae_2A$ and $C \cong A \otimes_k A$ as rings. 
Also, $C = Be_1B$ and $C \cong B \otimes_k B$ as rings.  
\end{proposition}
\begin{proof}
Clearly $Ae_2A \subseteq C$.  Conversely, if $c \in C$, then $c = \sum_j E_{M_1}(cu_j)v_j$.
But $\sum_j u_j \otimes v_j = \lambda^{-1}\sum_i z_ie_2 \otimes e_2 w_i$ by the endomorphism ring
theorem and the fact that both are dual bases for $E_{M_1}$.
Then
$
c = \lambda^{-1} \sum_i E_{M_1}(cz_i e_2)e_2 w_i \in Ae_2 A
$
as desired.

Since $e_2 w e_2 = E_M(w) e_2$ for every $w \in M_1$, we obtain the $E_M$-multiplication
on $Ae_2 A$.  Then $C = Ae_2 A = A \otimes_M A \cong A \otimes_k A$ since
$A \cap M = C_M(N) = k1_M$.

For the second statement, we observe:
\[  C = Ae_2A = A e_2 e_1 e_2 A \subseteq Ce_1C = B e_1 B,
\]
while the opposite inclusion is immediate.  
The ring isomorphism follows from the identity:
\begin{equation}
e_1 ce_1 = e_1 E_{M_1}(c)
\label{E_M_1 multiplication}
\end{equation}
for all $c \in C$, since $B \cap N1_2 \subseteq Z(N) = k1$.  
For there are $a_i,b_i \in A$ such that $c = \sum_i a_i e_2 b_i$ ,
and $\eta, \eta': A \rightarrow k$ such that, for all $a \in A$, 
 $e_1 a = e_1 \eta(a)$ while $ae_1 = \eta'(a)e_1$ by irreducibility.  Then 
we  easily compute  that $\eta = \eta'$.  Then:
\begin{eqnarray*}
e_1 c e_1 & = & \sum_i e_1 a_i e_2 b_i e_1 
           =  \sum_i \eta(a_i) \eta(b_i) e_1 e_2 e_1 \\
          & = & \lambda \sum_i e_1 a_i b_i = e_1 E_{M_1}(c).         \qed
\end{eqnarray*}
\renewcommand{\qed}{}\end{proof}

In Section~3 it will be apparent that $\eta$ is the counit 
$\eps$ on $A$. 

\begin{corollary}
\label{full matrix algebra}
If $n = \# \{ u_j \} = \# \{ v_j \} $, then 
$
C \cong M_n(k)
$
where the characteristic of $k$ does not divide $n$.
\end{corollary}
\begin{proof}
Since $B$ is a Frobenius algebra with Frobenius homomorphism
$E_{M_1}$, it follows from the isomorphism, $\End_k(B) \cong B \otimes B$ that  
\begin{equation}
C \cong \End_k(B) \cong M_n(k).
\end{equation}
We have $\mathrm{char}\, k \not | \, n$ since the index
$\lambda^{-1} = n1_k \neq 0$. 
\end{proof}

Since we can use $A$ in place of $B$ to conclude that
$C \cong \End_k(A)$ in the proof above, we see that $\dim_k A = \dim_k B$.
Although $C$ has a faithful trace, we will prefer 
the faithful linear functional $F$ studied below for its Markov-like properties in
Corollary~\ref{markov property of F}.
\begin{proposition}
\label{F}
$F := E_M \circ E_{M_1}$ is a faithful linear functional on $C$.
 \end{proposition}
\begin{proof}
We see that $E_M(E_{M_1}(C)) \in C_M(N) = k1$, and we identify $k1$ with $k$.
If $c = a \in A1_2$, we see that
\[
F(aC) = E_M(a E_{M_1}(C)) = E_M(aA) = 0
\] implies that $a = 0$ by Proposition~\ref{A is separable}, since
$E_M$ is a Frobenius homomorphism on $A$ and therefore faithful.

If $c = b \in B$, then by Lemma~\ref{C is A tensor B}
\[
F(bC) = E_M E_{M_1}(bC) = E_M(E_{M_1}(bB)A) = 0
\]
implies  first $E_{M_1}(bB) = 0$, hence $b = 0$. 

If $c \in C$, then there are $a_i \in A$ ($=E_{M_1}(cu_i)$) such that
$c = \sum_i a_i v_i$.  Then
\[
F(cC) = \sum_i E_M(a_i A)E_{M_1}(v_iB) = 0
\]
implies that each $a_i = 0$, since if $a_i \neq 0$, then $E_{M_1}(v_iB) = 0$, a contradiction.
Hence, $F$ is faithful on $C$.
\end{proof}
 
Denote the Nakayama automorphism of $F$ on $C$ by $q: C \rightarrow C$.  
It follows from Corollary~\ref{full matrix algebra} that $q$ is an inner automorphism.  
We note some other Nakayama automorphisms and study next their inter-relationships.
Let $q_A: A \rightarrow A$ be the Nakayama automorphism for $E_M$ on $A$.

Let $q_B: B \rightarrow B$ be the Nakayama automorphism for $E_{M_1}$ on $B$.
Let $\tilde{q}: B \rightarrow B$ be the Nakayama automorphism for 
$\hat{F} := E_M \circ E_{M_1}: M_2 \rightarrow M$,
a Frobenius homomorphism by Proposition~\ref{transitivity}. 

\begin{figure}
\[
\begin{array}{ccc}
 C & 
\xrightarrow{\scriptscriptstyle q} & C \\
 \downarrow  & &   \downarrow \\
  A & \stackrel{\scriptscriptstyle q_A}{\longrightarrow} & A  
\end{array}
\]
\caption{The  vertical arrows are given by the 
conditional expectation $E_{M_1}|_C$.}
\label{fig-em1}
\end{figure}

\begin{proposition}
\label{Nakayama autos}
We have $q_B = \tilde{q} = q|_B$, $q_A = q|_A$
and commutativity of the diagram in Figure~\ref{fig-em1}. 
\end{proposition}
\begin{proof}
We have for each $b \in B$, $c \in C$:
\[
F(cb) = F(q(b)c) = F(\tilde{q}(b) c)
\]
whence by faithfulness $q|_B = \tilde{q}$.  Then $q$ sends $B$ onto itself, so
\[
E_{M_1}(q_B(b_2)b_1)= E_{M_1}(b_1 b_2) = F(b_1b_2) = F(q(b_2)b_1) = E_{M_1}(q(b_2)b_1) 
\]
for each $b_1, b_2 \in B$, whence $q_B = q|_B$. 

 As for $q_A$, we note that
\[
F(q(a)c) = F(ca) = F(E_{M_1}(c)a) = F(q_A(a) E_{M_1}(c)) = F(q_A(a)c)
\]
for every $a \in A, c \in C$, whence $q = q_A$ on $A$. 

Commutativity of Figure~\ref{fig-em1} follows from the computation
applying Eq.\ (\ref{mark}): 
\[
F(q(E_{M_1}(c))c') = F(c'E_{M_1}(c)) = F(E_{M_1}(c')c) = F(q(c)E_{M_1}(c'))
= F(E_{M_1}(q(c))c'),  
\]
for all $c,c' \in C$.
\end{proof}

We now compute the 
conditional expectation of $C$ onto $B$, a lemma we will need in Section~3.

\begin{lemma}
\label{conditional expectation}
The map $E_B: C \rightarrow B$ defined by
$
E_B(c) = \sum_j F(cu_j) v_j
$
for all $c \in C$ is a conditional expectation. 
\end{lemma}
\begin{proof}
We first note that $E_B$ is the identity on $B$, since $E_{M_1}(bu_j) \in k1_1$, whence 
$
E_B(b) = \sum_j E_M(1_1) E_{M_1}(bu_j)v_j = b.
$
Since $E_M(E_{M_1}(cu_j)) \in k1$ for all $c\in C$, we have for each $b,b' \in B$:
\begin{eqnarray*}
E_B(be_1b') &=& \sum_j F(be_1 b'u_j) v_j 
=  \sum_j E_M(e_1 E_{M_1}(b'u_j q^{-1}( b)) v_j   \\
& = & \lambda \sum_j E_{M_1}(bb'u_j)v_j = \lambda bb'
\end{eqnarray*}
It follows from Proposition~\ref{C is A tensor A} that $E_B$ is a 
$B$-$B$-bimodule homomorphism (it corresponds to 
$\lambda \mu:  B \otimes B \rightarrow B$ under the  isomorphism 
$b \otimes b' \mapsto be_1b'$ of $B \otimes B$ with $C$).  

That $E_B$ is a Frobenius homomorphism follows 
from \cite[Lemma 2.6.1]{GHJ}, if we show it is one-sided faithful, e.g., $E_B(Cc) = 0$
implies $c = 0$. But this follows from $F$ being faithful and orthogonality of the dual 
bases  $\{ u_i \}$ and  $\{ v_i \}$.
\end{proof}
 
The corresponding conditional expectation $E_A: C \rightarrow A$ is easily seen
to be $E_{M_1}$ restricted to $C$.  We next record several Markov-like properties
of $F: C \rightarrow k$.

\begin{corollary}
\label{markov property of F}
The linear functional $F$ satisfies the following properties with respect to
$E_{M_1}$ and $E_B$:
\begin{eqnarray}
F(aE_{M_1}(c)) &=& F(ac),\qquad F(E_{M_1}(c)a) = F(ca), \label{mark} \\
F(bE_B(c)) &=& F(bc),\qquad  F(E_B(c)b) = F(cb),  \nonumber
\end{eqnarray}
for all $a\in A,b\in B,c\in C$. In particular, we have the following
Markov relations:
\begin{equation*}
F(ae_2) = F(e_2a) = \lambda F(a), \qquad F(be_1) =F(e_1b) =\lambda F(b).
\end{equation*}
\end{corollary}
\begin{proof}
According to the definitions of $F$ and $E_B$, we have
$F\circ E_{M_1} = F\circ E_B =  F$ and also $E_{M_1}(e_2) =\lambda$,
$E_B(e_1)=\lambda$, whence the result.
\end{proof}

\subsection*{The Pimsner-Popa identities.}  We note that:
\begin{eqnarray*}
\lambda^{-1} e_1 E_{M}(e_1x) &=&  e_1 x \ \ \ \forall x \in M_1 \\
\lambda^{-1} e_2 E_{M_1}(e_2 y) &=& e_2 y \ \ \ \forall y \in M_2.
\end{eqnarray*}
\textit{Proof}. Let  $x = \sum_i m_i \otimes m'_i$
where $m_i, m'_i \in M_1$.  Then $e_2 x = e_2 \sum_i E_M(m_i)m'_i$,
and $E_{M_1}(e_2x) = \lambda \sum_i E_M(m_i)m'_i$ from which one of the equations follows.
The other equation is similarly shown, as are
 the opposite Pimsner-Popa identities. 

\begin{corollary}
\label{q(e)}
$e_1\in Z(A)$, $e_2\in Z(B)$, and we have $q(e_1)=e_1$, $q(e_2)=e_2$.
\end{corollary}
\begin{proof}
From Eq.\ (\ref{E_M_1 multiplication}) 
$$
e_1a =  e_1 a e_1 = 
 ae_1,
$$
for all $a\in A$. It is clear from Eq.\ (\ref{Naka})
that a Nakayama automorphism fixes elements in the center
of a Frobenius \textit{algebra}. The assertions
about $e_2$ are shown similarly.
\end{proof}

\subsection*{When Hopf-Galois extensions are strongly separable.}

We recall a few facts about Hopf-Galois extensions \cite{M}. 
If $H$ is a finite dimensional
Hopf $k$-algebra with counit $\eps$
and comultiplication $\Delta(h) = h\1 \otimes h\2$,  then 
its dual $H^*$ is a Hopf algebra as well (and $H^{**} \cong H$).  Thus 
we have the following dual notions of algebra extension:    
$M/N$ is a right $H^*$-comodule algebra extension with coaction 
$M \rightarrow M \otimes H$, denoted by $\rho(a) = a\0 \otimes
a\1$, and $N = \{ b \in M |\, \rho(b) = b \otimes 1 \}$ if and only if
$M/N$ is a left $H$-module algebra extension with action of $H$ on $M$ given
by $h \lact a = a\0 \la a\1,\,h\ra$ 
and $N = \{  b \in M|\, \forall h \in H,\, h \lact b = \eps(h) b \}$.
Conversely, given an action of $H$ on $M$ and dual bases $\{ u_j\}$, $\{ p_j \}$
for $H$ and $H^*$, a coaction is given by 
\begin{equation}
\rho(a) = \sum_j (u_j \lact  a) \otimes p_j.
\label{coacting}
\end{equation}  

Recall on the one hand that 
$M/N$ is an $H^*$-\textit{Galois extension} if it is a right $H^*$-comodule
algebra such that the \textit{Galois map} $\beta: M \otimes _N M \rightarrow M \otimes H^*$
given by $a \otimes a' \mapsto a a'\0 \otimes a'\1$ is bijective.   

Recall
on the other hand that given a left $H$-module algebra $M$, there is the
smash product $M \# H$ with subalgebras $M = M \# 1$, $H = 1 \# H$ and 
commutation relation 
$ha = (h\1 \lact a)h\2$ for all $a \in M, h \in H$.  
If $N$ again denotes the subalgebra of invariants, then there
is a natural algebra homomorphism of the smash product into the right endomorphism
ring, $\Psi: M \# H \rightarrow \End(M_N)$ given 
by $m \# h \mapsto m( h \lact \cdot)$. 
We will  use the following basic proposition in Section~5 (and prove part
of the forward
implication below): 

\begin{proposition}[\cite{KT, U}]
\label{Ulbrich-Kreimer-Takeuchi}
An $H$-module algebra extension $M/N$ is $H^*$-Galois if and only if 
 $M \# H\xrightarrow{\cong} \End(M_N) $ via $\Psi$, and $M_N$ is a 
finitely generated projective module.
\end{proposition} 

The following theorem is  a converse
to our main theorem in~\ref{Galois extension}.
Let $H$ be a  finite dimensional,  semisimple and cosemisimple Hopf algebra. 

\begin{theorem}[Cf.\ \cite{K2}, 3.2]
\label{converse}
Suppose $M$ is a $k$-algebra  
and left $H$-module algebra
with subalgebra of invariants $N$.  If $M/N$ is an irreducible
 right $H^*$-Galois extension, then $M/N$ is a
strongly separable, irreducible extension of depth 2 with $\End(M_N) \cong M \# H$. 
\end{theorem}
\begin{figure}
\[
\begin{array}{ccc}
 M \otimes_N M & 
\xrightarrow{\scriptscriptstyle \beta} & M \otimes H^* \\
 \downarrow  & &   \downarrow {\scriptscriptstyle \cong}\\
  \End M_N & \xleftarrow{\scriptscriptstyle \Psi} & M \# H 
\end{array}
\]
\caption{Commutative diagram where the left vertical mapping is given
by $m \otimes m' \mapsto \lambda_m E \lambda_{m'}$ and the right vertical
mapping is the isomorphism $\id \otimes \theta$.}
\label{fig-beta}
\end{figure}
\begin{proof}
Since $H$ is finite dimensional (co)semisimple, $H$ is (co)unimodular and
there are integrals $f \in \int_{H^*}$ and $t \in \int_{H}$ such that
$f(t) = f(S(t)) = 1_k$, $\eps(t) = 1$ and $f(1) \neq 0$.  Moreover,
$g \mapsto (t \actl g)$ gives a Frobenius isomorphism  $\theta: H^* \xrightarrow{\cong}H $,
where $t \actl f = f(t\1)t\2 = 1_H$, since $f$ integral in $H^*$ means $x \actl f = f(x) 1_H$
for every $x \in H$.  

If $\beta:  M \otimes_N M \rightarrow M \otimes H^*$ is the Galois isomorphism,
given by $m \otimes m' \mapsto mm'\0 \otimes m'\1$, then $
\psi = (\id_M \otimes \theta) \circ
\beta$ is the isomorphism  $M \otimes_N M \xrightarrow{\cong} M \# H$  given by
\begin{eqnarray*} 
 \ m \otimes m' \mapsto mm'\0 \otimes (t \actl m'\1) 
&=& m\la m'\1 , t\1 \ra m'\0 \otimes t\2 \\
&=& m(t\1 \cdot m') \otimes t\2 = mtm'.
\end{eqnarray*}

Now define $E: M \rightarrow N$ by $E(m) = t \cdot m$, where $t \cdot m \in N$ since
$h\cdot (t \cdot m) = (ht)\cdot m = \eps (h) t \cdot m$.  Note
that $E$ is an $N$-$N$-bimodule map and $E(1) = \eps(t)1 = 1$.

Denote $\beta^{-1}(1 \otimes f) = \sum_i x_i \otimes y_i \in M \otimes_N M$.
Since $(\id \otimes \theta)(1 \otimes f) = 1 \# 1$, which is sent by $\Psi$ to
$\id_M$, it follows that $\sum_i x_i (Ey_i) = \id_M$ (cf.\ Figure~\ref{fig-beta}).\footnote{
  $E$ is in fact an Frobenius homomorphism
with dual bases $\{ x_i \}$, $\{ y_i \}$, the other equation,  
$\sum_i (x_i E) y_i = \id_M$, following readily from a computation
using  $\beta' = \eta \circ \beta$, where $\beta'$ is the ``opposite'' Galois
mapping given by $\beta'(m \otimes m') = m\0 m'\otimes m\1$ and
$\eta$ is an automorphism of $M \otimes H^*$ given by $\eta(m \otimes g) = m\0 \otimes m\1 S(g)$
 \cite{KT}.}

The homomorphism $\Psi: M \# H \rightarrow \End(M_N)$ (given by $
m \# h \longmapsto (m' \mapsto m(h \cdot m'))
$) is now readily checked to have 
 inverse mapping given by $g \mapsto
\sum_i g(x_i) t y_i$ \cite{KT}. 

By counitarity of the $H^*$-comodule $M$, then $\mu: M \otimes_N M
\rightarrow M$ factors through $\beta$ and the map $M \otimes H^* \rightarrow M$
given by $m \otimes g \mapsto mg(1)$.  Then $\sum_i x_i  y_i = f(1_H) 1_M$,
whence the $k$-index $[M:N]_E$ is $\lambda^{-1} = f(1_H)$.  

It is not hard to compute that $C_{M \# H}(N) = C_M(N) \# H$ which is $H$ since
$M/N$ is irreducible. Since $M \# H$ is free over $M$ with basis in $H$, 
we see that the first half of the depth 2 condition is satisfied. 

The second half of depth 2 follows from noting that $M \# H$ is a right $H$-Galois extension
of $M$. For the coaction $ M \# H \rightarrow (M \# H) \otimes H$
is given by 
\begin{equation}
m \# h \mapsto m \# h\1 \otimes h\2.
\label{eq:endo}
\end{equation}
One may compute the inverse of the Galois map to be given by
 $\beta^{-1}(m \# h \otimes h') = mhS(h'\1) \otimes h'\2$. Then 
$
M_2 \cong M \# H \# H^*  $ and the rest of the proof proceeds as in the previous paragraph. 
\end{proof}

The proof shows that an $H$-Galois extension $M/N$ has an endomorphism ring
theorem: ${\mathcal E}/M$ is an $H^*$-Galois extension. 
A converse to the endomorphism ring theorem depends on  ${\mathcal E}/M$
being \textit{$H^*$-cleft},
as discussed in Section~6. 

\end{section}

%%%%%%%%%%%%%%%%%%%%%%%%%%%%%%%%%%%%%%%%%%%%%%%%%%%%%%%%%%%%%%%%%%%%%%%%%%%%%

\begin{section}
{Hopf algebra structures on centralizers}

In this section, we define and study an important non-degenerate pairing
of $A$ and $B$ given by Eq.\ (\ref{pairing}).  This transfers the algebra
structure of $A$ onto a coalgebra structure of $B$, and vice versa.
The rest of the section is devoted to showing that $B$ is a Hopf algebra
with an antipode $S$ satisfying $S^2 = \id$. The key step in this 
and the next sections is Proposition~\ref{exchange relation}.  
 
\subsection*{A duality form}

As in Section~2, we let 
$N\subset M \subset M_1 \subset M_2 \subset \cdots$
be the Jones tower constructed from a strongly separable
irreducible extension $N\subset M$ of depth $2$,
$F =E_M\circ E_{M_1}$ denote the functional on $C$ 
defined in Proposition~\ref{F},
$e_1  \in M_1,\, e_2\in M_2$
be the first two Jones idempotents of the tower, and 
$\lambda^{-1} =[M:N]$ be the index.

\begin{proposition}
%[Cf.\ \cite{S}, Proposition 10]
\label{duality}
The  bilinear form
\begin{equation}
\la a, \, b \ra = \lambda^{-2} F(ae_2e_1b), \qquad a\in A,\, b\in B,
\label{pairing}
\end{equation}
is non-degenerate on $A\otimes B$.
\end{proposition}
\begin{proof}
If $ \la a,\, B\ra =0$ for some $a\in A$, then
we have $F(ae_2e_1c)=0$ for all $c\in C$,
since $e_1B = e_1 C$ by Lemma~\ref{e_2A is e_2C}.
Taking $c = e_2q^{-1}(a') \,(a'\in A)$ and using 
the braid-like relations between
Jones idempotents and Markov property 
(Corollary~\ref{markov property of F}) of $F$ we have 
$$
F(a'a) = \lambda^{-1}F(a'ae_2) = \lambda^{-1}F(ae_2q^{-1}(a')) =
\lambda^{-2}F(ae_2e_1(e_2 q^{-1}(a')) = 0
$$
for all  $a'\in A$, therefore $a=0$ (by Proposition~\ref{A is  separable}). 

Similarly, if $\la A,\, b\ra =0$ for some $b$, then 
$F(ce_2e_1b)=0$ for all $c\in C$, which for $c= q(b')e_1\,(b'\in B)$ gives
$$
F(bb') = \lambda^{-1}F(e_1bb') = \lambda^{-1}F(q(b')e_1 b) =
\lambda^{-2}F( (q(b')e_1) e_2e_1 b) = 0
$$
for all $b'\in B$, therefore  $b=0$. 
\end{proof}

Observe that since $k$ is a field the  Proposition above 
shows that the map $b\mapsto E_{M_1}(e_2e_1b)$ is
a linear isomorphism between $B$ and $A$.
Indeed, $E_{M_1}(e_2e_1b) = 0$ 
implies that for all $a \in A$ one has 
\[
F(ae_2e_1b) = F(aE_{M_1}(e_2e_1b)) = 0, 
\]
whence $b=0$ by nondegeneracy.
 
\subsection*{A coalgebra structure}
Using the above duality form we introduce a coalgebra structure on $B$.

\begin{definition}
\label{coalgebra}
The algebra $B$ has a comultiplication 
$\Delta: B \rightarrow B \otimes B$,  $b \mapsto b\1 \otimes b\2$ defined by
\begin{equation}
\la a,\, b\1 \ra \la a',\, b\2 \ra = \la a a',\, b\ra
\end{equation} 
for all $a,a'\in A,\, b\in B$, and counit $\eps: B \rightarrow k$
given by ($\forall b \in B$)
\begin{equation}
\eps(b) = \la 1,\,b \ra.
\end{equation} 
\end{definition}

\begin{proposition}
\label{1 and eps}
For all $b,c\in B$ we have :
\begin{equation}
\eps(b) = \lambda^{-1}F(be_2),
\end{equation}
\begin{equation}
\Delta(1) = 1\otimes 1,
\end{equation}
\begin{equation}
\eps(bb') =\eps(b)\eps(b').
\end{equation}
\end{proposition}
\begin{proof}

We use the Pimsner-Popa identities 
together with Corollaries~\ref{markov property of F} and \ref{q(e)}
to compute
\begin{eqnarray*}
\eps(b) 
&=& \lambda^{-2} F(e_2e_1b) = 
\lambda^{-2} F(e_1be_2) = \lambda^{-1}F(be_2),\\
\la a, 1\ra \la a', 1\ra
&=& \lambda^{-4}F(ae_2e_1)F(a'e_2e_1) \\
&=& \lambda^{-2} F(ae_1) F(a'e_1) = \lambda^{-2} F(aE_M(a'e_1)e_1) \\
&=& \lambda^{-1} F(aa'e_1) = \la aa',\, 1\ra,\\
\eps(b)\eps(b') 
&=& \lambda^{-2} F(be_2)F(b'e_2) = \lambda^{-2} F(b E_{M_1}(b'e_2)e_2) \\
&=& \lambda^{-1} F(bb'e_2) = \eps(bb'),
\end{eqnarray*}
for all $a,a'\in A,\, b,b'\in B$
(note that the restriction of $E_M|_A = F$ and
$E_{M_1}|_B =F$, identifying $k$ and $k1$).
\end{proof}

\subsection*{The antipode of $B$}

Recall that the map $b\mapsto E_{M_1}(e_2e_1b)$ is
a linear isomorphism between $B$ and $A$. But
considering the Jones tower $N^{op}\subset M^{op} \subset M_1^{op} 
\subset M_2^{op}$ of the opposite algebras, 
we conclude that the map $b\mapsto E_{M_1}(be_1e_2)$ is a linear
isomorphism as well. This lets us  define a linear map 
$S: B\to B$, called the \textit{antipode}, as follows.

\begin{definition}
\label{antipode}
For every $b\in B$ define $S(b)\in B$ to be the unique element
such that
$$
F(q(b)e_1e_2a) = F(ae_2e_1S(b)), \qquad \mbox{ for all } a\in A,
$$
or, equivalently,
$$
E_{M_1}(be_1e_2) = E_{M_1}(e_2e_1S(b)). 
$$

\end{definition}

\begin{remark}
\label{more on S}
Note that $S$ is bijective and that the above condition  implies 
\begin{equation}
E_{M_1}(bxe_2) = E_{M_1}(e_2xS(b)), \qquad \mbox { for all } x\in M_1.
\label{remarkeq}
\end{equation}
Indeed, $B$ commutes with $M$ and any $x\in M_1$ can be written  
as $x = \sum_i\, m_ie_1n_i$ with $m_i,n_i \in M$, so that
$$
E_{M_1}(bxe_2) = \Sigma_i\, m_i E_{M_1}(be_1e_2) n_i 
= \Sigma_i\, m_i E_{M_1}(e_2e_1S(b))  n_i = E_{M_1}(e_2xS(b)).
$$
\end{remark}

\subsection*{$A$ and $B$ are Hopf algebras}
To prove that $B$ is Hopf algebra, it remains to show that
$\Delta$ is a homomorphism and that $S$ satisfies the antipode axioms.
The next proposition is also the key ingredient 
for an action of $B$ on $M_1$ which makes
$M_2$  a smash product.

\begin{proposition}
\label{exchange relation}
For all $b\in B$ and $y\in M_1$ we have
$$
yb =  \lambda^{-1} b\2 E_{M_1}(e_2yb\1).
$$
\end{proposition}
\begin{proof}
First, let us show that the above equality holds true in
the special case $y = e_1$. Let $E_B$ be the  
conditional expectation from $C$ to $B$
given by $E_B(c) = \Sigma_i\, F(cu_i)v_i$ as 
in Proposition~\ref{conditional expectation}.
  
We claim that for any $c\in C$
we have $c=0$ if  $\la a,\, E_B(ca')\ra =0$
for all $a,a'\in A$. For since $C = BA$, let $c=\sum_i\, b_i a_i$ with
$a_i\in A$ and $b_i\in B$, then
$$
\la a,\, E_B(ca')\ra = \sum_i\, \la a,\,b_i E_B(a_ia')\ra
= \sum_i\, \la a,\,b_i \ra F(a_ia'),
$$
and the latter expression is equal to $0$ for all $a,a'\in A$
only if for each $i$ either $a_i=0$ or $b_i=0$.

Observe that $q$ restricted to $A$
coincides with the Nakayama automorphism $q_A : A\to A$
of the Frobenius extension $M_1/N$ since
$$
F(q(a)c) = F(ca) = F(E_{M_1}(c)a) =E \circ E_M(q_A(a)E_{M_1}(c)) = F(q_A(a)c),
$$
therefore, using the Pimsner-Popa identity for $C = Be_1B$,  
we establish the proposition for $y = e_1$: 
\begin{eqnarray*}
\la a,\, E_B(e_1 b a')\ra
&=& \lambda^{-2} F(ae_2e_1 E_B(e_1 b a')) \\
&=& \lambda^{-1}  F(ae_2e_1 ba') =  \lambda \la q(a')a,\, b\ra,\\
\la a,\, \lambda^{-1} b\2 E_B(E_{M_1}(e_2e_1b\1)a') \ra
&=& \lambda^{-1} \la a,\, b\2  \ra F(e_2e_1b\1a') \\
&=& \lambda \la a,\, b\2  \ra \la q(a'),\, b\1 \ra =
\lambda \la q(a')a,\, b\ra,
\end{eqnarray*}
since $E_B|_A = F$. 

Next, arguing as in Remark~\ref{more on S} we write
$y =\Sigma_i\, m_ie_1n_i$ with $m_i,n_i \in M$, whence
$$
yb = \Sigma_i\, m_ie_1bn_i = \lambda^{-1} \Sigma_i\,
m_i b\2 E_{M_1}(e_2e_1b\1) n_i = b\2 E_{M_1}(e_2yb\1). \qed 
$$
\renewcommand{\qed}{}\end{proof}

\begin{corollary}
\label{action}
For all $b\in B$ and $x,y \in M_1$ we have:
$$
E_{M_1}(e_2xyb) = \lambda^{-1} E_{M_1}(e_2xb\2)E_{M_1}(e_2yb\1).
$$
\end{corollary}
\begin{proof}
The result follows from multiplying the identity
from Proposition~\ref{exchange relation} by $e_2x$
on the left and taking $E_{M_1}$ from both sides.
\end{proof}

Although the antipode axiom
(cf.\ Prop.~\ref{S is an antipode}) implies 
that $S$ is a coalgebra anti-homomorphism, we will have
to establish these two properties of  $S$ in the reverse order, 
as stepping stones to Propositions~\ref{Delta is a homo} 
and \ref{S is an antipode}.

\begin{lemma}
\label{anti-coalgebra}
$S$ is a coalgebra anti-automorphism.
\end{lemma}
\begin{proof}
For all $a,a'\in A$ and $b\in B$ we have by Corollary~\ref{action} :
\begin{eqnarray*}
\la aa',\, S(b) \ra
&=& \lambda^{-2} F(q(b)e_1e_2aa') 
= \lambda^{-3} F(e_1e_2 E_{M_1}(e_2aa'b)) \\
&=& \lambda^{-4} F(e_1e_2 E_{M_1}(e_2ab\2) E_{M_1}(e_2a'b\1) )\\
&=& \lambda^{-6} F(e_1e_2 E_{M_1}(e_2ab\2) ) F(e_1e_2  E_{M_1}(e_2a'b\1) ) \\
&=& \lambda^{-4} F(e_1e_2ab\2) F(e_1e_2 a'b\1) \\
&=& \lambda^{-4} F(q(b\2)e_1e_2a) F(q(b\1)e_1e_2a') \\
&=& \la a, S(b\2) \ra \la a', S(b\1) \ra, 
\end{eqnarray*}
where we use  the definition of $S$, the Pimsner-Popa identity, and
Corollary~\ref{markov property of F}. 
 Thus, $\Delta(S(b))=S(b\2) \otimes S(b\1)$.
\end{proof}

\begin{corollary}
\label{left action}
For all $b\in B$ and $x,y\in M_1$ we have :
$$
E_{M_1}(bxye_2) = \lambda^{-1} E_{M_1}(b\1xe_2)E_{M_1}(b\2ye_2)
$$
\end{corollary}
\begin{proof}
We obtain this formula by replacing $b$ with $S(b)$ in 
Corollary~\ref{action}
and using Eq.\ (\ref{remarkeq}) as well as Lemma~\ref{anti-coalgebra}.
\end{proof}

\begin{proposition}
\label{S squared}
$S^2 = q|_B^{-1}$.
\end{proposition}
\begin{proof}
The statement follows from the direct computation :
\begin{eqnarray*}
F(ae_2e_1q^{-1}(b))
&=& \lambda^{-1} F(E_{M_1}(bae_2)e_2e_1) \\
&=& \lambda^{-1} F(E_{M_1}(e_2a S(b))e_2e_1) \\
&=& \lambda^{-1} F(e_2 E_{M_1}(e_2a S(b))e_1) \\
&=& F(e_2aS(b) e_1) = F( a E_{M_1}( S(b)e_1e_2) ) \\
&=& F(ae_2e_1S^2(b)),
\end{eqnarray*}
for all $a\in A$ and $b\in B$,
using Remark~\ref{more on S} and Corollary~\ref{q(e)}.
\end{proof}

\begin{proposition}
\label{Delta is a homo}
$\Delta$ is an algebra homomorphism.
\end{proposition}
\begin{proof}
Note that $q|_B$ is a coalgebra automorphism by Proposition~\ref{S squared}.

By Corollary~\ref{left action} we have, for all $a,a'\in A$ 
and $b,b'\in B$ :
\begin{eqnarray*}
\la aa',\, bb' \ra 
&=& \la \lambda^{-1}E_{M_1}(q(b')aa'e_2),\, b \ra \\
&=& \la \lambda^{-2}E_{M_1}(q(b')\1 ae_2)E_{M_1}(q(b')\2 a'e_2),\, b \ra \\
&=& \la \lambda^{-1}E_{M_1}(q(b'\1) ae_2), b\1 \ra
    \la \lambda^{-1}E_{M_1}(q(b'\2)  a' e_2),\, b\2 \ra    \\
&=& \la a,\, b\1 b'\1 \ra \la a',\, b\2 b'\2 \ra,
\end{eqnarray*}
 whence $\Delta(bb') =  \Delta(b) \Delta(b')$.
\end{proof}

\begin{proposition}
\label{S is an antipode}
For all $b\in B$ we have $S(b\1)b\2 =\eps(b)1 = b\1 S(b\2)$.
\end{proposition}
\begin{proof}
Using Corollary~\ref{left action} and the definition of the antipode
we have 
\begin{eqnarray*}
\la a,\, S(b\1) b\2 \ra
&=& \lambda^{-1} \la E_{M_1}(q(b\2) ae_2),\,S(b\1) \ra \\
&=& \lambda^{-3} F(q(b\1) e_1 e_2 E_{M_1}(q(b\2) ae_2)) \\
&=& \lambda^{-3} F( E_{M_1}(q(b\1) e_1 e_2 ) E_{M_1}(q(b\2) ae_2)) \\
&=& \lambda^{-2} F(q(b) e_1 a  e_2) = \lambda^{-2} F(e_1 a  e_2 b) \\
&=& \lambda^{-2} F(e_1 a) F(b e_2) = \la a,\, 1\eps(b) \ra ,
\end{eqnarray*}
$\forall a\in A,\, b\in B$. The second identity follows
similarly from  Corollary~\ref{action} 
and the corollary $q \circ S = S^{-1}$ from Proposition~\ref{S squared}:
\begin{eqnarray*} 
\la a,\,b\1 S(b\2) \ra
&=& \lambda^{-1} \la E_{M_1}( q(S(b\2))ae_2),\, b\1 \ra \\
&=& \lambda^{-3} F( E_{M_1}( q(S(b\2))ae_2) e_2 e_1 b\1) \\
&=& \lambda^{-3} F( E_{M_1}(S^{-1}(b\2)ae_2) e_2 e_1 b\1) \\ 
&=& \lambda^{-3}F( E_{M_1}(e_2 a b\2) E_{M_1}(e_2 e_1 b\1)) \\
&=& \lambda^{-2} F(e_2 a  e_1 b) = \lambda^{-2} F(ae_1 be_2) = \la a,\, 1\eps(b) \ra,
\end{eqnarray*}
i.e., $S$ satisfies the antipode properties.
\end{proof}

\begin{theorem}
\label{Hopf algebra}
$A$ and $B$ are semisimple  
Hopf algebras.
\end{theorem}
\begin{proof}
Follows from Propositions~\ref{1 and eps}, \ref{Delta is a homo},
\ref{S is an antipode}, and \ref{A is separable}.
Note that semisimplicity and separability are notions that coincide
for finite dimensional Hopf algebras \cite{M}. 
The non-degenerate duality form of Proposition~\ref{duality}
makes $A$ the Hopf algebra dual to $B$.    
\end{proof}

\begin{corollary}
The antipodes of $A$ and $B$ satisfy $S^2 = \id$, $F$ is a trace,
and $A$, $B$ are Kanzaki strongly separable.
\label{EG consequence}
\end{corollary}
\begin{proof}
Etingof and Gelaki  proved that a semisimple and cosemisimple
Hopf algebra is involutive \cite{EG}.  
It follows from Proposition~\ref{S squared} that $q_B = \id_B$. But we compute:
$$
\la a,\, q^{-1}(b)\ra = \lambda^{-2}F(bae_2e_1) = \lambda^{-2}F(q(a)e_2e_1b)
= \la q(a),\,b\ra
$$
for all $a \in A$, $b \in B$, from which it follows that $S_A^2 = q_A = \id_A$.
Since $C = BA$, we have  $q = \id_C$. Whence 
$F$, $E_M$ and $E_{M_1}$ are traces on $C$, $A$ and $B$, respectively.

It follows from Proposition~\ref{A is separable} that $\{ \lambda^{-1} E_{M_1}|_B, \lambda u_i, v_i \}$
is a separable basis for $B$; similarly, $\{ \lambda^{-1} E_M|_A, \lambda z_i, w_i \}$
is a separable basis for $A$, whence $A$ and $B$ are strongly separable
algebras. 
\end{proof}  

\begin{remark}
Note that $e_2$ is a ($2$-sided) integral in $B$, since
$
\la a, e_2 b \ra = \la a, e_2 \ra \eps(b) = \la a, be_2 \ra
$
by the Pimsner-Popa identity.  Similarly, $e_1$ is an integral in $A$.   
\end{remark}

\end{section}

%%%%%%%%%%%%%%%%%%%%%%%%%%%%%%%%%%%%%%%%%%%%%%%%%%%%%%%%%%%%%%%%%%%%%%%%%%

\begin{section}
{Action of $B$ on $M_1$ and $M_2$ as a smash product}

In this section we define the Ocneanu-Szyma\'{n}ski action of $B$ on $M_1$, which
makes $M_1$ a $B$-module algebra (cf.\ Eq.\ (\ref{action B on M_1})).  We then
 describe $M$ as its subalgebra of invariants and
$M_2$ as the smash product algebra of $B$ and $M_1$.
As a corollary, we note that $M_1/M$ and $M_2/M_1$ are respectively
$A$- and $B$-Galois extensions.  

\begin{proposition}
%[Cf.\ \cite{S}, Proposition 17]
\label{action of B}
The map $\lact : B \otimes M_1 \to M_1$ :
\begin{equation}
b \lact x = \lambda^{-1} E_{M_1}(b x e_2)
\label{action B on M_1}
\end{equation}
defines a left $B$-module algebra action on $M_1$,
called the \textit{Ocneanu-Szyma\'{n}ski action}.
\end{proposition}   
\begin{proof}
The above map defines a left $B$-module structure
on $M_1$, since $1\lact x = \lambda^{-1} E_{M_1}(x e_2) = x$ and
$$
b \lact (c \lact x)
= \lambda^{-2}  E_{M_1}(b  E_{M_1}(c x e_2)  e_2)
= \lambda^{-1}  E_{M_1}( bc x e_2 ) = (bc) \lact x.
$$
Next, Corollary~\ref{left action} implies that  
$b\lact xy = ( b\1 \lact x)( b\2 \lact y)$.
Finally, $b\lact 1 = \lambda^{-1} E_{M_1}(b e_2) = 
 \lambda^{-1} F(b e_2)1 = \eps(b)1$. 
\end{proof} 

We note that an application of Proposition~\ref{exchange relation} provides
another formula for the action of $B$ on $M_1$: 
\begin{equation}
b \lact x = b\1 x S(b\2).
\label{outeraction} 
\end{equation} 

%The coaction $\rho: M_1 \rightarrow M_1 \otimes B^*$ is given by
%$ \rho(m) = \sum_j (u_j \lact m) \otimes E_{M_1}v_j|_B$. 
%

\begin{proposition}
\label{invariants}
$M_1^B =M$, i.e., $M$ is the subalgebra of invariants of $M_1$. 
\end{proposition}
\begin{proof}
If $x\in M_1$ is such that $b\lact x = \eps(b)x$ for all
$b \in B$, then $ E_{M_1}(bxe_2) = \lambda\eps(b) x$. Letting 
$b=e_2$ we obtain $ E_{M}(x) =  \lambda^{-1} E_{M_1}(e_2 x e_2)
=  \eps(e_2) x =x$, therefore $x \in M$.

Conversely, if $x\in M$, then $x$ commutes with $e_2$  and
$$
b\lact x = \lambda^{-1} E_{M_1}(b e_2 x) =
\lambda^{-1} E_{M_1}(be_2)x = \eps(b)x,
$$
therefore $M_1^B = M$. 
\end{proof}

Note from the proof that  $ e_2 \lact x = E_M(x)$, i.e.,
the conditional expectation $E_M$ is action on $M_1$ by the 
integral $e_2$ in $B$. The rest of this section is strictly speaking
not required for Section~6. 

\begin{theorem}
\label{smash product}
The map $\theta : x\# b \mapsto xb$ defines an algebra  isomorphism
between the smash product algebra $M_1\# B$ and $M_2$. 
\end{theorem}
\begin{proof}
The bijectivity of $\theta$ follows from Lemma~\ref{tensor product}.
To see that $\theta$ is a homomorphism it suffices to 
note that $ by = (b\1\lact y) b\2$ for all $b\in B$ and $y\in M_1$.
Indeed, using Eq.\ (\ref{outeraction}), 
\begin{eqnarray*}
(b\1\lact y) b\2
&=& b\1 y S(b\2)b\3 \\
&= & b\1 y \eps(b\2) = by. \qed
%&=& \lambda^{-1} E_{M_1}(b\1ye_2)b\2 \\
%&=& \lambda^{-2} b\3 E_{M_1}(e_2E_{M_1}(b\1ye_2)b\2) \\
%&=& \lambda^{-1} b\3 E_{M_1}(e_2 y S(b\1)b\2) = by.  \qed
\end{eqnarray*}
\renewcommand{\qed}{}\end{proof}

From this and Lemma~\ref{C is A tensor B}, we conclude that:
\begin{corollary}
\label{C is A smash B}
$C \cong A \# B.$
\end{corollary}

\begin{corollary}
\label{Galois extension}
$M_1/M$ is an $A$-Galois extension.  $M_2/M_1$ is a $B$-Galois extension. 
\end{corollary}
\begin{proof}
Dual to the  left $B$-module algebra $M_1$ defined above
 is a right $A$-comodule algebra
$M_1$ with the same subalgebra of coinvariants $M$,
since $B^* \cong A$. 
By Theorem~\ref{smash product}  and the endomorphism ring theorem, 
$M_1 \# B \xrightarrow{\cong} M_2 \xrightarrow{\cong}
\End^r_M(M_1)$ is given by the natural map $x \# b  \mapsto x( b \lact \cdot)$
since if $b = \sum_i a_i e_2 a'_i$ for $a_i,a'_i \in A$, then for all 
$y \in M_1$,  
$$
x(b \lact y) =  \lambda^{-1} \sum_i x a_i E_{M_1}(e_2 a'_i y e_2) 
              = x \sum_i a_i E_M(a'_i y) .
$$
By Proposition~\ref{Ulbrich-Kreimer-Takeuchi} then, 
$M_1$ is a right  $A$-Galois extension of $M$. 

It follows from the endomorphism ring theorem for Hopf-Galois extensions
(cf.\ end of Section~3) that $M_2/M_1$ is $B$-Galois.   
\end{proof} 

Since $M_2$ is a smash product of $M_1$ and $B$, thus a $B$-comodule
algebra, it has
a left $A$-module algebra action given by applying Eq.\ (\ref{eq:endo}):
$$  a \lact (mb) = \la a, b\2 \ra mb\1,$$
for every $a \in A, m \in M_1, b \in B$.\footnote{Alternatively, the depth
2 condition is satisfied by $M_1/M$ due to Theorem~\ref{converse},
and $C_{M_3}(M_1) \cong A$ via $a \mapsto d$ where $F(ae_2e_1b) = 
E_{M_1}E_{M_2}(be_3e_2d)$
for all $b \in B$; 
whence we may repeat the arguments in Sections 3 -- 5 to define an $A$-module
algebra action on $M_2$, 
$a \lact m_2 = \lambda^{-1} E_{M_2}(dm_2e_3)$, where $M_3$, $E_{M_2}$
and $e_3$ are of course  the basic construction of $M_2/M_1$.
This is the same action of $A$ on $M_2$ by  repeating Proposition~\ref{total
integral}.} We remark that 
$M_1/M$ and $M_2/M_1$ are   \textit{faithfully flat} (indeed free)  
Hopf-Galois extensions with normal basis property \cite{M}[chap.\ 8].

\end{section}
%%%%%%%%%%%%%%%%%%%%%%%%%%%%%%%%%%%%%%%%%%%%%%%%%%%%%%%%%%%%%%%%%%%%%%%%%
\begin{section}
{Action of $A$ on $M$ and $M_1$ as a smash product}

In this section, we note that $M_1/M$ is an $A$-cleft $A$-extension
(Proposition~\ref{total integral}).  It follows from a theorem in the Hopf
algebra literature
that $M_1$ is a crossed product of $M$ and $A$.  The cocycle $\sigma$ determining
the algebra structure of $M \#_{\sigma} A$ is in this case trivial.  Whence
$M_1 \cong M \# A$ and $M/N$ is a left $B$-Galois extension
(Theorem~\ref{M_1 is a smash product}).  
We end the section with a proof
of Theorem~\ref{non-commutative analogue} and a proposal for
further study.  

From the Ocneanu-Szyma\'{n}ski action given in Eq.\ (\ref{action B on M_1}), we note that $B \lact A = A$. 
The next proposition  shows, based on Corollary~\ref{EG consequence}, that the action of $B$ on $A$
yields a coaction $ A \rightarrow A \otimes A$ (when dualized)
which is identical with the comultiplication on $A$. Recall that
an extension of $k$-algebras $N' \subseteq M'$ is called an \textit{$A$-extension} if
$A$ is a Hopf algebra co-acting on $M'$ such that $M'$ is a right $A$-comodule algebra
with $N' = {M'}^{{\rm co}A}$ \cite{M}:  e.g. $M_1/M$ is an $A$-extension by duality since $A$
is finite dimensional.  An $A$-extension $M'/N'$ is \textit{$A$-cleft} if there is a right
$A$-comodule map $\gamma: A \rightarrow M'$ which is invertible with respect to the convolution
product on $\Hom(A,M')$ \cite{M,DT}. 

\begin{proposition}
\label{total integral}
The natural inclusion $\iota: A \hookrightarrow
M_1$ is a total integral such that the $A$-extension $M_1/M$ is $A$-cleft.
\end{proposition}
\begin{proof}
Since $\iota(1) = 1$, we  show that $\iota$ is a total integral by showing it is 
a right $A$-comodule morphism \cite{DT}. 
Denoting the coaction $M_1 \rightarrow M_1 \otimes A$ (which
is the  dual of Action~\ref{action B on M_1}) by
$w \mapsto w\0 \otimes w\1$, we have  $w\0 \la w\1,b \ra = b \lact w$ 
for every $b \in B$.  Since each $a\0 \in A$ by Eq.\ (\ref{coacting}),
it  suffices to check that
$a\0 \otimes a\1 = a\1 \otimes a\2$:
\begin{eqnarray*}
\la a\1,b \ra \la a\2,b' \ra & = & \la a,bb' \ra 
     =  \lambda^{-2} F(ae_2 e_1 bb' ) \\
     & = & \lambda^{-3} F(E_{M_1}(b'ae_2)e_2 e_1 b)  
     =  \la \lambda^{-1} E_{M_1}(b'ae_2),b \ra  \\
     & = & \la a\0,b \ra \la a\1,b' \ra. 
\end{eqnarray*}
Finally, we note that $\iota$ has convolution inverse in $\Hom(A,M_1)$
given by $\iota \circ S$ where $S: A \rightarrow A$ denotes the antipode
on $A$. 
\end{proof}

We recall the following result of Doi and Takeuchi (see also \cite[Prop.\ 7.2.3]{M}
and \cite{BCM}):
\begin{proposition}[\cite{DT}]
\label{Doi-Takeuchi}
Suppose $M'/N'$ is an $A$-extension, which is $A$-cleft by a total integral 
$\gamma: A \rightarrow M'$. Then there is a crossed product
action of $A$ on $N'$ given by
\begin{equation}
a \cdot n = \gamma(a\1)n \gamma^{-1}(a\2)
\label{DT}
\end{equation}
for all $a \in A, n \in N'$, and a  cocycle $\sigma: A \otimes A \rightarrow N'$
given by
\begin{equation}
\sigma(a,a') = \gamma(a\1)\gamma(a'\1)\gamma^{-1}(a\2 a'\2)
\end{equation}
for all $a,a' \in A$, such that $M'$ is isomorphic as algebras to a crossed product of
$A$ with $N'$ and cocycle $\sigma$:  $$M' \cong N' \#_{\sigma} A$$
given by $n \# a \mapsto n \gamma(a)$.
\end{proposition}

Applied to our $A$-cleft $A$-extension $M_1/M$, we conclude:

\begin{theorem}
\label{M_1 is a smash product}
$M_1 $ is isomorphic to the smash product $M \# A$ via $m \# a \mapsto ma$.
\end{theorem}
\begin{proof}
The cocycle $\sigma$ associated to $\iota: A \rightarrow M_1$ is trivial,
since $$\sigma(a,a') = a\1 a'\1 S(a\2 a'\2) = \eps(a) \eps(a') 1_1.$$
It follows from Eq.\ (\ref{DT}) and  \cite[Lemma 7.1.2]{M})
that $M$ is an $A$-module algebra with action $A \otimes M \rightarrow M$ given by 
\begin{equation}
a \lact m = a\1 m S(a\2).
\label{action of A on M}
\end{equation}
It follows from Proposition~\ref{Doi-Takeuchi} and triviality of the crossed 
product that $M_1$ is a smash product of $M$ and $A$ as claimed. 
\end{proof}

\begin{lemma}
The fixed point algebra is $M^A = N$.
\label{fixpoint}
\end{lemma}
\begin{proof}
That $N \subseteq M^A$
follows from the definition of $A$ and its Hopf algebra structure. 
Conversely, suppose that 
$m \in M$ is such that $a \lact m = \eps(a) m$ for all $a \in A$.
In a computation similar to that of \cite{So}, 
we note that $am = ma$ in $M_1$ for any $a \in A$:
\[ 
am  =  a\1 m S(a\2)a\3 
    =  (a\1 \lact m) a\2 =  ma. \]
Letting $a = e_1$, we see that $m$ commutes with $e_1$, so that
$E(m) e_1 = e_1 m e_1 = e_1 m$.  Applying $E_M$ to this, we arrive
at $m = E(m) \in N$. 
\end{proof}

\begin{theorem}
$M/N$ is a $B$-Galois extension.
\label{M is B-Galois over N}
\end{theorem}
\begin{proof}
This follows from Theorem~\ref{M_1 is a smash product} 
and Proposition~\ref{Ulbrich-Kreimer-Takeuchi}, if we 
 prove that $\Psi: M \# A \rightarrow \End(M_N)$ given
by $$m \# a \mapsto (x \mapsto m (a \lact x))$$
is an isomorphism.

Towards this end, we claim that $e_1 \lact x = E(x)$ for every $x \in M$.
Let $G = e_1 \lact \cdot$.  A few short calculations using
 Lemma~\ref{fixpoint} show 
 that $G \in \Hom_{N-N}(M,N)$ such that $G|_N = \id_N$, since 
$$ae_1 = \lambda^{-1}E_M(ae_1)e_1 = \lambda^{-2}F(ae_2e_1)e_1 = \eps_A(a)e_1$$
and 
$$\eps_A(e_1)
= \lambda^{-2}F(e_1e_2e_1) = 1. $$ Since $E$ freely generates
$\Hom_N(M,N)$ (as a Frobenius homomorphism), there is $d \in C_M(N) = k1$
such that $G = Ed$, whence $E = G$ as claimed.

Then $\Psi((m \# e_1)( m' \# 1_A)) = \lambda_m E \lambda_{m'}$ for 
all $m, m' \in M$ is
 surjective. An inverse mapping may be defined by $f \mapsto
\sum_i (f(x_i) \# e_1)(y_i \# 1_A)$ for each
$f \in \End(M_N)$, where $\{ x_i \}$,
$\{ y_i \}$ are dual bases for $E$ as in Section~2. 
\end{proof} 

We are now in a position to  note  the proof of 
Theorem~\ref{non-commutative analogue}. 

\begin{theorem}[= Theorem~\ref{non-commutative analogue}]
If $M/N$ is an irreducible  extension of depth 2, then $M/N$
is  strongly separable if and only if $M/N$ is an $H$-Galois extension, where
$H$ is a semisimple, cosemisimple Hopf algebra. 
\end{theorem}
\begin{proof}
The forward implication  follows from Theorem~\ref{M is B-Galois over N}.   
The reverse implication follows  from Theorem~\ref{converse}. 
\end{proof}

We propose the following two problems related to this paper:
\begin{enumerate}
\item  Are conditions 1 and 2 in the depth 2 conditions independent?
\item  What is a suitable definition of normality for $M/N$ extending the notion of normal
 field extensions?\footnote{There is an early definition of non-commutative normality in 
\cite{E}. }
\end{enumerate}

\end{section}

%%%%%%%%%%%%%%%%%%%%%%%%%%%%%%%%%%%%%%%%%%%%%%%%%%%%%%%%%%%%%%%%%%%%%%%%%%

\bibliographystyle{amsalpha}

\end{document}